\documentclass[10pt,oneside]{article}
\usepackage[OT4]{fontenc}\usepackage[cp1250]{inputenc}\usepackage{a4} \usepackage{amsfonts}\usepackage{amsmath}\usepackage{amssymb}\usepackage{amsxtra}
\usepackage{amsthm}

\usepackage[dvips]{graphicx}\usepackage{pstricks}\usepackage{color}
\usepackage{bezier}\usepackage{epic}\usepackage{eepic}
\usepackage{pstricks}
\usepackage{pst-node}
\usepackage{lineno}
\usepackage[normalem]{ulem}\usepackage{epsf}

\newtheorem{thm}{Theorem}
\newtheorem{cor}{Corollary}
\newtheorem{df}{Definition}

\newtheorem{obs}{Observation}

\newtheorem{lem}{Lemma}

\theoremstyle{remark}

\newcommand{\cp}{\, \square \,}

\let\oldenumerate\enumerate
\renewcommand{\enumerate}{
  \oldenumerate
  \setlength{\itemsep}{0.5pt}
  \setlength{\parskip}{0pt}
  \setlength{\parsep}{0pt}
}

\newcommand{\1}{ \vspace{0.1cm} }
\newcommand{\2}{ \vspace{0.2cm} }

\begin{document}

\title{Graphs with disjoint 2-dominating sets}
\author{$^{1}$Michael A. Henning\thanks{Research
supported in part by the University of Johannesburg} \, and \, $^{2}$Jerzy Topp\\
\\
$^1$Department of Mathematics and Applied Mathematics\\
University of Johannesburg \\
Auckland Park 2006, South Africa \\
\small {\tt Email: mahenning@uj.ac.za} \\
\\
$^2$State University of Applied Sciences\\ 82-300 Elbl\c{a}g, Poland \\
\small {\tt Email: j.topp@pwsz.elblag.pl}}

\date{}
\maketitle

\begin{abstract}
A subset $D \subseteq V(G)$ is a dominating set of a multigraph $G$ if every vertex in $V(G) \setminus D$ has a neighbor in $D$, while $D$ is a 2-dominating set of $G$ if every vertex belonging to $V(G)\setminus D$ is joined by at least two edges with a vertex or vertices in $D$. A graph $G$ is a $(2,2)$-dominated graph if it has a pair $(D,D')$ of disjoint $2$-dominating sets of vertices of $G$. In this paper we present two characterizations of minimal $(2,2)$-dominated graphs.
\end{abstract}

{\small \textbf{Keywords:} Domination; 2-domination} \\
\indent {\small \textbf{AMS subject classification: 05C69, 05C85}}

\section{Introduction}

For notation and graph theory terminology we generally follow~\cite{Henning-Yeo2013}. Let $G = (V(G),E(G))$ be a~graph with possible multi-edges and multi-loops, and with vertex set $V(G)$ and edge set $E(G)$. For a vertex $v$ of $G$, its {\em neighborhood\/}, denoted by $N_{G}(v)$, is the set of vertices adjacent to $v$. The {\em closed neighborhood\/} of $v$, denoted by $N_{G}[v]$, is the set $N_{G}(v)\cup \{v\}$.  In general, for a subset $X\subseteq V(G)$, the {\em neighborhood\/} of $X$, denoted by $N_{G}(X)$, is defined to be $\bigcup_{v\in X}N_{G}(v)$, and the {\em closed\/} neighborhood of $X$, denoted by $N_{G}[X]$, is the set $N_{G}(X)\cup X$. The $2$-\emph{neighborhood} of $v$, denoted by $N_G^2(v)$, is the set of vertices at distance~$2$ from $v$ in $G$, that is, $N_G^2(v) = \{u \in V(G) \colon d_G(u,v) = 2\}$. The \emph{closed} $2$-\emph{neighborhood} of $v$, denoted by $N_G^2[v]$, is the set of vertices within distance~$2$ from $v$ in $G$, and so $N_G^2[v] = N_G[v] \cup N_G^2(v)$.

If $A$ and $B$ are disjoint sets of vertices of $G$, then we denote by $E_G(A,B)$ the set of edges in $G$ joining a vertex in $A$ with a~vertex in $B$. For one-element sets we write $E_G(v,B)$,  $E_G(A,u)$, and $E_G(u,v)$ instead of $E_G(\{v\},B)$, $E_G(A,\{u\})$, and $E_G(\{u\},\{v\})$, respectively. If $v$ is a~vertex of $G$, then by $E_G(v)$ we denote the set of edges incident with $v$ in $G$. The {\em degree} of a vertex $v$ in $G$, denoted by $d_G(v)$, is the number of edges incident with $v$ plus twice the number of loops incident with $v$. A vertex of degree one is called a {\em leaf}. A vertex is {\em isolated\/} if its degree equals zero. For an integer $k \ge 1$, we let $[k] = \{1,\ldots,k\}$.

A set of vertices $D \subseteq V(G)$ of $G$ is a {\em dominating set\/} if every vertex in $V(G)\setminus D$  has a~neighbor in $D$, while $D$ is a {\em $k$-dominating set\/}, where $k$ is a positive integer,  if every vertex belonging to $V(G)\setminus D$ is joined by at least $k$ edges with a vertex or vertices in $D$. If $G$ is a graph without multiple edges, then a~subset $D\subseteq V(G)$ is a $k$-{\em dominating set} of $G$ if $|N_G(v)\cap D|\ge k$ for every $v\in V(G)\setminus D$.

If $k$ and $\ell$ are positive integers, then a pair $(D_1,D_2)$ of proper and disjoint subsets of the vertex set $V(G)$ of a graph $G$ is a {\em $(k,\ell)$-pair\/} in $G$ if $D_1$ is a $k$-dominating set of $G$, and $D_2$ is an $\ell$-dominating set of $G$. A graph $G$ is said to be a {\em $(k,\ell)$-dominated graph\/} if it contains a $(k,\ell)$-pair. It is obvious from the above definition, that if a graph $G$ is a $(k,\ell)$-dominated graph, then necessarily $\max\{k,\ell\}\le \Delta(G)$, $1\le \min\{k,\ell\}\le \delta(G)$, and $k+\ell\le |V(G)|$. Trivially, if $G$ is a $(k,\ell)$-dominated graph, then $G$ is a $(k',\ell')$-dominated graph, where $1\le k'\le k$ and $1\le \ell' \le \ell$. In addition, if $G$ is a~$(k,\ell)$-dominated graph, then $G$ is an $(\ell,k)$-dominated graph. Thus we may suppose that if $G$ is~a $(k,\ell)$-dominated graph, then $k\le \ell$.

We observe that a complete graph $K_n$ is a~$(k,\ell)$-dominated graph (for positive integers $k$ and $\ell$) if and only if $k+\ell\le n$. Moreover, we observe that a complete bipartite graph $K_{m,n}$ is a $(m,n)$-dominated graph. A cycle $C_n$ is a $(2,2)$-dominated graph if and only if $n$ is an even positive integer, while every cycle of odd length is a $(1,2)$-dominated graph but not a $(2,2)$-dominated graph.

\begin{figure}[h!]
\begin{center} \bigskip
{\epsfxsize=5in \epsffile{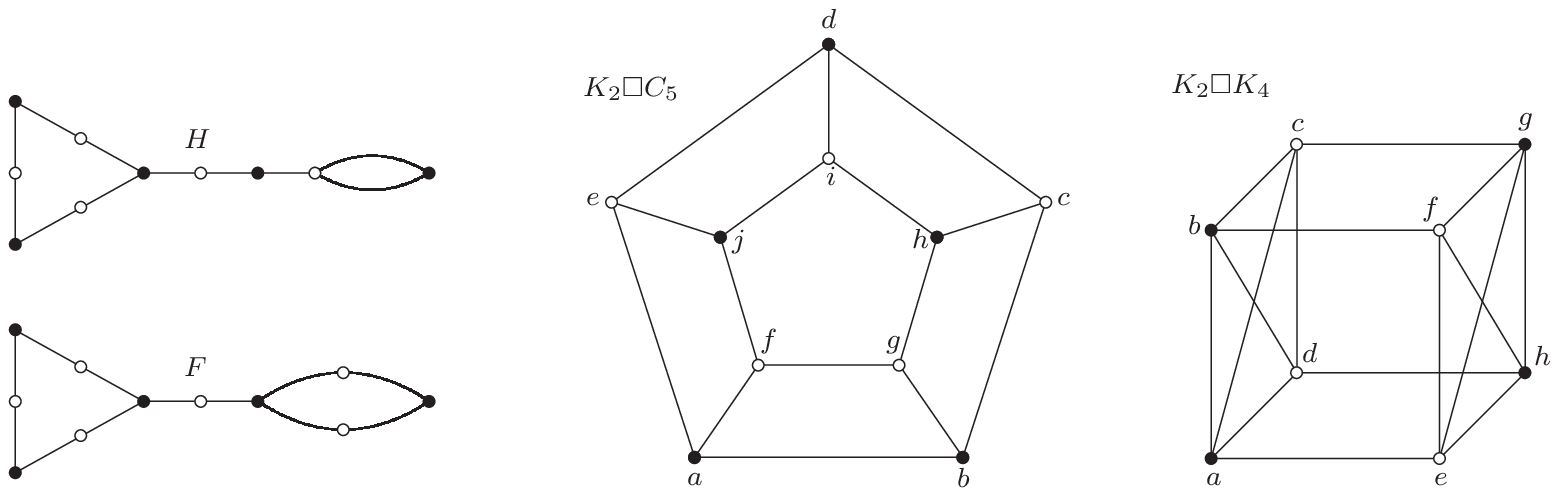}}

\vspace{-2mm}
\caption{Graphs $F$, $H$, $K_2\Box C_5$, and $K_2\Box K_4$}\label{K2naK4iC5} \label{pierwsze-grafy}
\end{center}
\end{figure}

Of the graphs in Fig.~\ref{K2naK4iC5}, the graphs $F$, $H$, and the Cartesian product $K_2 \cp C_5$ are examples of $(2,2)$-dominated graphs, while the Cartesian product $K_2 \cp K_4$ is an example of a $(3,3)$-dominated graph (and a $(1,4)$-dominated graph). The appropriate $(2,2)$- and $(3,3)$-pairs in these graphs are determined by the sets of black and white vertices, respectively, illustrated in Fig.~\ref{K2naK4iC5}. More generally, we show in Corollary~\ref{cor:Cartesian-product} that if $G$ and $H$ are graphs without isolated vertices, then the Cartesian product $G \cp H$ is a $(2,2)$-dominated graph.

\vspace{-6mm}

Ore \cite{Ore1962} was the first to observe that a graph without isolated vertices contains two disjoint dominating sets. That is, Ore observed that every such graph is a $(1,1)$-dominated graph. Subsequently, various properties of graphs having disjoint dominating sets of different types have been extensively studied, for example, in papers \cite{Bang-Jensen Bessy Havet Yeo}--\cite{Miotk-Zylinski},  to mention just a few. All $(1,2)$-dominated graphs were characterized in
\cite{Henning-Rall2013,Miotk-Topp-Zylinski,Miotk-Zylinski}. In this paper, we study $(2,2)$-dominated graphs, and, in particular,  we present two characterizations of minimal $(2,2)$-dominated graphs. It is worth mentioning here that it follows from \cite[Theorem 12]{Bang-Jensen Bessy Havet Yeo} that in the general case it is $\cal NP$-complete to decide whether a~given graph $G$ is a~$(2,2)$-dominated graph if $\delta(G)\ge 2$.

\section{Elementary properties of $(2,2)$-dominated graphs}

In this section, we present properties of $(2,2)$-dominated graphs that will need in order to prove our main results.

\begin{df}
{\rm A connected graph $G$ is said to be a \emph{minimal $(2,2)$}-\emph{dominated graph}, if $G$ is a~$(2,2)$-dominated graph and no proper spanning subgraph of $G$ is a $(2,2)$-dominated graph. }
\end{df}

From this definition, we immediately have the following observations.

\begin{obs}
\label{observ-2-supergraph}
Every spanning supergraph of a $(2,2)$-dominated graph is a $(2,2)$-dominated graph, and every $(2,2)$-dominated graph is a spanning supergraph of some~minimal $(2,2)$-dominated graph.
\end{obs}

\begin{obs}
\label{observ-1-bipartite graph}
Every bipartite graph of degree at least~$2$ $($and every spanning supergraph of such a graph$)$ is a $(2,2)$-dominated graph.
\end{obs}

As a consequence of Observation~\ref{observ-1-bipartite graph}, we have the following result.

\begin{cor}
\label{cor:Cartesian-product}
If $G$ and $H$ are graphs without isolated vertices, then their Cartesian product $G \cp H$ is a $(2,2)$-dominated graph.
\end{cor}
\begin{proof} We may assume that $G$ and $H$ are connected graphs each of order at least~$2$. Let $T_G$ and $T_H$ be spanning trees of $G$ and $H$, respectively. Then $T_G \cp  T_H$ is a bipartite spanning subgraph of $G \cp H$ and $\delta(G \cp H) = 2$. Thus, by Observation \ref{observ-1-bipartite graph}, $T_G\cp T_H$ and $G\cp H$ are $(2,2)$-dominated graphs.
\end{proof}

In view of Observation~\ref{observ-2-supergraph}, minimal $(2,2)$-dominated graphs can be viewed as skeletons of $(2,2)$-dominated graphs, skeletons which can be extended to any $(2,2)$-dominated spanning supergraph.

The next theorem presents general properties of minimal $(2,2)$-dominated graphs.

\begin{thm}
\label{thm-property-22}
A graph $G$ is a minimal $(2,2)$-dominated graph if and only if $G$ has the following three properties. \\ [-24pt]
\begin{enumerate}
\item $\delta(G)\ge 2$.
\item $G$ is a bipartite graph.
\item Every edge of $G$ is incident with a vertex of degree~$2$ in $G$.
\end{enumerate}
\end{thm}
\begin{proof} Assume first that $G$ is a minimal $(2,2)$-dominated graph, and let $(D_1,D_2)$ be a $(2,2)$-pair in $G$. Since $D_1$ and $D_2$ are disjoint, every vertex $v$ of $G$ belongs to $V(G) \setminus D_1$ or to $V(G)\setminus D_2$, and thus $|N_G(v)\cap D_1|\ge 2$ or $|N_G(v)\cap D_2|\ge 2$ (since $D_1$ and $D_2$ are 2-dominating sets), implying that $\delta(G)\ge 2$. We now claim that $D_1$ and $D_2$ form a partition of $V(G)$. Suppose, to the contrary, that $V(G) \setminus (D_1\cup D_2)\ne \emptyset$. Then, for every $v\in V(G)\setminus (D_1\cup D_2)$, the pair $(D_1\cup \{v\},D_2)$ is a~$(2,2)$-pair in $G-E_G(v,D_1)$ (and $(D_1,D_2\cup \{v\})$ is a~$(2,2)$-pair in $G-E_G(v,D_2)$), a contradiction to the minimality of $G$. From the minimality of $G$ it also follows that $G$ is a bipartite graph in which the sets $D_1$ and $D_2$ form a bipartition, for if two vertices $x$ and $y$ belonging to $D_1$ (or $D_2$) were adjacent in $G$, then $(D_1,D_2)$ would be a $(2,2)$-pair in $G-xy$, a contradiction. Finally, no two vertices of degree at least~$3$ are adjacent in $G$, for if a vertex $x\in D_1$ of degree at least~$3$ were adjacent to a vertex $y\in D_2$ of degree at least~$3$, then $(D_1,D_2)$ would be a $(2,2)$-pair in $G-xy$.  From this and from the fact that $\delta(G)\ge 2$, it follows that every edge of $G$ is incident with a vertex of degree 2 in $G$ (and, therefore, $\delta(G)=2$).

Assume now that $G$ is a bipartite graph with partite sets $A$ and $B$, in which $\delta(G) \ge 2$ and every edge of $G$ is incident with a vertex of degree~$2$ in $G$. Then  $(A,B)$ is a $(2,2)$-pair in $G$ and therefore $G$ is a $(2,2)$-dominated graph. Now, if $e$ is an edge of $G$, then $G-e$ has a vertex of degree 1 (since $e$ is incident with a vertex of degree 2) and therefore $G-e$ is not a $(2,2)$-dominated graph. Consequently, $G$ is a minimal $(2,2)$-dominated graph. \end{proof}

If $H$ is a graph (with possible multi-edges or multi-loops), then the \emph{subdivision graph} of $H$, denoted by $S(H)$, is the graph obtained from $H$ by inserting a new vertex into each edge and each loop of $H$. We remark that the graphs $F$ in Fig. \ref{pierwsze-grafy}, $G$ in Fig. \ref{pierwszeprzykladsciagania}, and $G$ in Fig. \ref{graphs in F} are examples of subdivision graphs. We note that the subdivision graph $S(H)$ of $H$ is a bipartite graph. On the other hand, we have the following useful observation.

\begin{obs}
A connected graph $G$ is a subdivision graph if and only if $G$ is a connected bipartite graph with partite sets $A$ and $B$ such that at least one of them consists only of vertices of degree~$2$. Furthermore, a connected bipartite graph $G$ with partite sets $A$ and $B$ such that $\delta(G)\ge 2$ and $|B| \ge |A|$ is a subdivision graph if and only if $d_G(x)=2$ for every $x\in B$.
\end{obs}

We state the next two important corollaries of Theorem~\ref{thm-property-22} that will prove very helpful to us. This corollary states that every minimal $(2,2)$-dominated graph (and therefore every $(2,2)$-dominated graph) can be constructed from a subdivision graph.

\begin{cor}
If a minimal $(2,2)$-dominated graph has multi-edges, then at least one of the vertices incident with them is of degree~$2$. \end{cor}

\begin{cor}
\label{S(H) is minimal}
If $H$ is a graph with $\delta(H)\ge 2$ and with possible multi-edges or multi-loops, then its subdivision graph $S(H)$ is a minimal $(2,2)$-dominated graph.
\end{cor}

\section{Constructive characterization of minimal $(2,2)$-dominated graphs}%

We remark that both graphs $F$ and $H$ in Fig. \ref{pierwsze-grafy} are minimal $(2,2)$-dominated graphs. But only $F$ is a subdivision graph. Thus, not every minimal $(2,2)$-dominated graph is a subdivision graph.  Surprisingly, there are interesting connections between minimal $(2,2)$-dominated graphs and subdivision graphs. To prepare the ground for our explanation, let us begin with the following definition of a $\cal P$-contraction, which will play an important role in our considerations.

Let $G$ be a bipartite graph. We define a vertex $v$ as a \emph{contractible vertex}
of $G$ if $v$ is not incident with a multi-edge. Let $v$ be a contractible vertex
in $G$, and let ${\cal P}(v)$ be a partition of the neighborhood $N_G(v)$ of $v$.
Recall that $N_G^2(v)$ is the set of vertices at distance~$2$ from $v$ in $G$,
while $N_G^2[v]$ is the set of vertices within distance~$2$ from $v$ in $G$.
Let $G'=G({\cal P}(v))$ denote a graph in which
\[
V(G')= (V(G) \setminus N_G(v)) \cup (\{v\}\times {\cal P}(v)),
\]
and where
\begin{align*}
N_{G'}(u)  =\,\, & N_G(u) \hspace*{0.25cm} \mbox{if} \hspace*{0.25cm} u \in V(G') \setminus N_G^2[v], \\[1ex]
N_{G'}\big((v,S)\big)  =\,\, & N_G(S) \hspace*{0.25cm} \mbox{if} \hspace*{0.25cm} (v,S) \in \{v\} \times {\cal P}(v), \\[1ex]
N_{G'}(v) = \,\, & \{v\} \times {\cal P}(v) \hspace*{0.25cm} \mbox{and} \hspace*{0.25cm} |E_{G'}(v,(v,S))|=1  \hspace*{0.25cm} \mbox{for each}  \hspace*{0.25cm} S \in {\cal P}(v), \\[-2ex]
\intertext{and}\\[-6ex]
N_{G'}(u) =\,\,& \{(v,S) \colon S \in {\cal P}(v) \hspace*{0.15cm} \mbox{and} \hspace*{0.15cm}  N_G(u) \cap S \ne \emptyset\} \cup (N_G(u) \setminus N_G(v)) \1
\end{align*}
for every vertex $u \in N_G^2(v)$. Moreover, in this case when $u \in N_G^2(v)$ and $(v,S) \in N_{G'}(u)$, then $|E_{G'}(u,(v,S))| = |E_{G}(u,S)|$.

\begin{figure}[h!]
\begin{center} \bigskip
{\epsfxsize=5in \epsffile{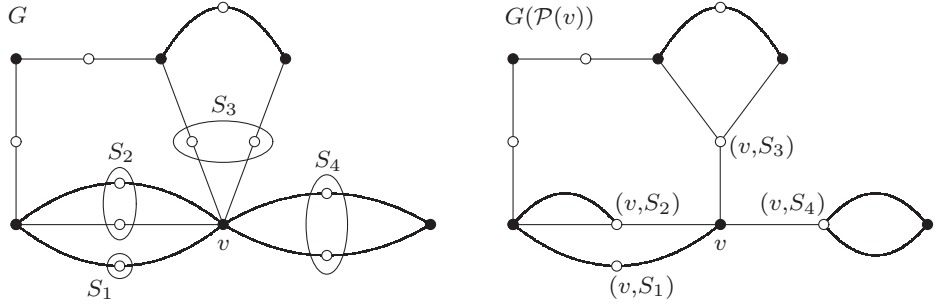}}

\vspace{-2mm}
\caption{Graph $G$ and its ${\cal P}$-contraction $G({\cal P}(v))$} \label{pierwszeprzykladsciagania}
\end{center}
\end{figure}

The graph $G({\cal P}(v))$ is called a ${\cal P}$-{\em contraction} of $G$ with respect to the partition ${\cal P}(v)$.  To illustrate this construction, we present on the left side of Fig.~\ref{pierwszeprzykladsciagania} a graph $G$ with a specified vertex $v$ and a~partition ${\cal P}(v)=\{S_1,S_2,S_3,S_4\}$ of the neighborhood $N_G(v)$ of $v$ into four subsets indicated by ellipses. The graph on the right side of Fig.~\ref{pierwszeprzykladsciagania} is the associated ${\cal P}$-contraction $G({\cal P}(v))$ of $G$ with respect to the partition ${\cal P}(v)$.

\vspace{-1mm}

The following observation follows readily from the definition of the ${\cal P}$-contraction of a graph.

\begin{obs}
\label{observation-on-P-contraction}
If $G$ is a bipartite graph and ${\cal P}(v) =\{S_1, \ldots, S_k\}$ is a partition of the neighborhood $N_G(v)$ of a contractible vertex $v$ of $G$, then the following properties hold in the ${\cal P}$-contraction $G' = G({\cal P}(v))$ of $G$. \\ [-22pt]
\begin{enumerate}
\item $G'$ is a bipartite graph.
\item $d_{G'}(v)= |{\cal P}(v)|=k$.
\item $d_{G'}((v,S_i))= 1+\sum_{u\in S_i}\big(d_G(u)-1\big)$ for every $S_i\in {\cal P}(v)$.
\item $d_{G'}(x) = d_G(x)$ for every $x \in V(G') \setminus N_G^2[v]$.
\item $G'$ is isomorphic to $G$ if $|{\cal P}(v)|=d_G(v)$, that is, if ${\cal P}(v)$ consists of singletons.
\end{enumerate}
\end{obs}

We are interested in determining partitions ${\cal P}(v)$ of $N_G(v)$ which transform
a minimal $(2,2)$-dominated graph $G$ into a minimal $(2,2)$-dominated graph $G({\cal P}(v))$. We begin with the following lemma.

\begin{lem} \label{partition of N(v)}
Let $G$ be a minimal $(2,2)$-dominated graph, and let ${\cal P}(v)$ be a partition of $N_G(v)$ for some contractible vertex $v$ of $G$, say ${\cal P}(v)= \{S_1,\ldots, S_k\}$, where $1\le |S_1|\le \cdots \le |S_k|$.  Then the ${\cal P}$-contraction $G({\cal P}(v))$ of $G$ is a minimal $(2,2)$-dominated graph if and only if at least one of
the following two statements holds.
\\ [-22pt]
\begin{enumerate}
\item $k = |N_G(v)|$.
\item $k = 2$ and $d_G(x)=2$ for every $x\in N_G(S_i)\setminus \{v\}$ if $|S_i|\ge 2$ $(i\in \{1,2\})$.
\end{enumerate}
\end{lem}

\begin{proof} From the fact that $G$ is a minimal $(2,2)$-dominated graph and from Theorem \ref{thm-property-22} it follows that $G$ is a bipartite graph, $\delta(G)=2$, and every edge of $G$ is incident with a vertex of degree~$2$. Let $G'$ denote the ${\cal P}$-contraction $G({\cal P}(v))$ of $G$, where $v$ is a contractible vertex of $G$, ${\cal P}(v)= \{S_1,\ldots, S_k\}$ is a partition of $N_G(v)$ and $1\le |S_1|\le \cdots \le |S_k|$.

We shall show that $G'$ is a minimal $(2,2)$-dominated graph if and only if at least one of the statements (a) and (b) holds. Since the result is obvious if $k = |N_G(v)|$
(as in this case ${\cal P}(v)= \{x\colon x\in N_G(v)\}$ and $G'$ is isomorphic to $G$, see Observation \ref{observation-on-P-contraction}(e)), we may assume that $k <|N_G(v)|$. Then $k<|N_G(v)|=|S_1|+\ldots+|S_k|\le k|S_k|$ and therefore $|S_k|\ge 2$. In addition, it follows from Theorem~\ref{thm-property-22} that neither the case $k=1<|N_G(v)|$ nor the case
$3\le k<|N_G(v)|$ is possible as otherwise either $v$ is of degree 1 in $G'$ or $v$
and $(v,S_k)$ are adjacent vertices of degree at least three in $G'$. Thus it remains to consider the case $k=2<|N_G(v)|$.

It follows from Observation \ref{observation-on-P-contraction}\,(a)--(d) that $G'$ is a
bipartite graph and $\delta(G')=2$, because $G$ is bipartite, $\delta(G)=2$, and $k=2$.
Consequently, by Theorem~\ref{thm-property-22}, to prove that $G'$ is a minimal $(2,2)$-dominated graph, it suffices to show that every edge of $G'$ is incident with a vertex of degree~$2$. Since the edges $v(v,S_1)$ and $v(v,S_2)$ are incident with $v$, which is of degree~$2$ in $G'$, and every edge of $G'$, which is not incident with $(v,S_1)$ or $(v,S_2)$, has inherited this property from the graph $G$, the graph $G'$ is a minimal $(2,2)$-dominated graph if and only if every edge of $G'$  incident with $(v,S_1)$ or $(v,S_2)$ (and different from $v(v,S_1)$ and $v(v,S_2)$) is incident a vertex of degree~$2$ in $G'$ (and in $G$). This property holds if and only if $|S_1|=1$ and $d_G(x)=2$ for every $x\in N_G(S_2)\setminus \{v\}$ or $2\le |S_1|\le |S_2|$ and $d_G(x)=2$ for every $x\in N_G(S_1\cup S_2)\setminus \{v\}$, that is, if and only if $d_G(x)=2$ for every $x\in N_G(S_i)\setminus \{v\}$ if $|S_i|\ge 2$ where $i\in \{1,2\}$. This completes the proof. \end{proof}

\medskip
By ${\cal M}$ we denote the family of all connected minimal $(2,2)$-dominated graphs. We are now in a~position to present a constructive characterization of the family ${\cal M}$. For this purpose, let ${\cal F}$ be the family  of graphs that:
\begin{itemize}
\item[$(1)$] contains the subdivision graph $S(H)$ for every connected graph $H$ with $\delta(H)\ge 2$ (and possibly with multi-edges and multi-loops); and
 \item[$(2)$] is closed under ${\cal P}$-contractions, that is, if a graph $G$ belongs to $\cal F$, then the ${\cal P}(v)$-contraction $G({\cal P}(v))$ of $G$ belongs to $\cal F$, if ${\cal P}(v)=\{S_1,S_2\}$ is a partition of $N_G(v)$, where $v$ is a contractible vertex of degree at least~$3$ in $G$ and where $|S_1| = 1$, and every vertex belonging to $N_G(S_2) \setminus \{v\}$ is of degree~$2$.
 \end{itemize}

\begin{figure}[h!]
\begin{center} \bigskip
{\epsfxsize=5in \epsffile{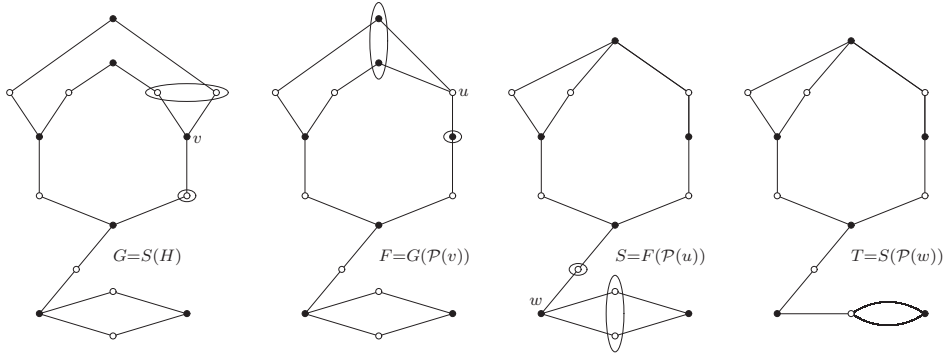}}

\vspace{-2mm}
\caption{Graphs $G$, $F$, $S$, and $T$ belonging to the family ${\cal F}$}
\label{graphs in F}
\end{center}
\end{figure}

Examples of graphs $G=S(H)$, and ${\cal P}$-contractions $F=G({\cal P}(v))$, $S=F({\cal P}(u))$, and $T=~S({\cal P}(w))$ belonging to the family ${\cal F}$ are given in Fig. \ref{graphs in F}.

The following theorem provides a characterization of minimal $(2,2)$-dominated graphs.

\begin{thm}
\label{char:M}
A connected graph $G$ is in the family $\cal M$ if and only if $G$ is in the family $\cal F$.
\end{thm}
\begin{proof} It follows from Corollary~\ref{S(H) is minimal} and Lemma~\ref{partition of N(v)} that ${\cal F} \subseteq {\cal M}$. Thus it remains to prove that ${\cal M} \subseteq {\cal F}$. Assume that $G$ is a connected graph belonging to ${\cal M}$. By Theorem \ref{thm-property-22}, $G$ is a~bipartite graph with $\delta(G)=2$ such that every edge of $G$ is incident with a vertex of degree~$2$ in $G$. Let $A$ and $B$ be the partite sets of $G$. Let

\begin{align*}
A_G^2 = \,\,& \{x\in A\colon d_G(x)=2\}\,\,  \mbox{and}\,\, 
A_G^3 = A\setminus A_G^2, \2 \\[-2ex]
\intertext{and}\\[-7ex]
B_G^2 = \,\,& \{x\in B\colon d_G(x)=2\}\,\,  \mbox{and}\,\, 
B_G^3 = B\setminus B_G^2.
\end{align*}

By induction on $k=\min\{|A_G^3|, |B_G^3|\}$ we will prove that $G\in {\cal F}$. If $k=0$, then at least one of the sets $A$ and $B$ consists of vertices of degree~$2$, implying that $G$ is a subdivision graph and proving that $G$ belongs to ${\cal F}$. Thus, let $k$ be a~positive integer, and assume that $|A_G^3| \ge |B_G^3|=k$.

Among all vertices $u \in A_G^3$ and $v\in B_G^3$, let $u$ and $v$ be chosen to be at minimum distance apart in $G$, that is, $d_G(u,v)= \min\{d_G(x,y) \colon x \in A_G^3\,\,\mbox{and}\,\, y\in B_G^3\}$. Let $P \colon v = v_0, v_1, \ldots, v_\ell = u$ be a shortest $(u,v)$-path in $G$. Since $u$ and $v$ belong to different partite sets of $G$, we note that $\ell$ is odd. Further since $G$ is a minimal $(2,2)$-dominated graph, the set $A_G^3 \cup B_G^3$ is an independent set, implying that $\ell \ge 3$. By the choice of the path $P$, every internal vertex of the path $P$ has degree~$2$ in $G$, while the vertices $u$ and $v$ are both of degree at least~$3$. Further from the minimality of $G$ and by Theorem~\ref{thm-property-22}, every neighbor of $v$ has degree~$2$ in $G$.

Without loss of generality we assume that the subset $N_G(v) \setminus \{v_1\}$ of $N_G(v)$ is the union of two disjoint sets $\{p_1,\ldots,p_m\}$ and $\{s_1, \ldots, s_n\}$, where each vertex $p_i$ has degree~$2$ in $G$ and is joined by a pair of parallel edges with $v$ (say by edges $e_i$ and $e_i^{'}$), while each vertex $s_j$ has degree~$2$ in $G$ and is adjacent to a vertex, say $s_j'$, different from $v$. We remark that possibly $s_i' = s_j'$ if $i \ne j$, and possibly one of the sets $\{p_1,\ldots,p_m\}$ and $\{s_1,\ldots, s_n\}$ is empty. Now let $G'$ be a graph with vertex set $V(G') = (V(G)\setminus \{v\})\cup V^*$, where
\[
V^* = \{(v,e_1), (v,e_1'), \ldots, (v,e_m), (v,e_m')\} \cup \{(v,s_1),\ldots, (v,s_n)\},
\]
and with edge set $E(G')$ obtained from $E(G)$ as follows:
\\ [-22pt]
\begin{enumerate}
\item[$\bullet$] deleting all edges incident with~$v$ in $G$,
\item[$\bullet$] adding an edge from $v_1$ to every vertex in the set $V^*$,
\item[$\bullet$] adding an edge from $p_i$ to both the vertices $(v,e_i)$ and $(v,e_i')$  for all $i \in [m]$, and
\item[$\bullet$] adding an edge from $s_i$ to the vertex $(v,s_i)$ for all $i \in [n]$.
\end{enumerate}
That is, defining
\[
\begin{array}{lcl}
E_1 & = & E(G) \setminus E_G(v), \1 \\
E_2 & = & \{v_1x \colon x \in V^* \}, \1 \\
E_3 & = & \{(v,e_i)p_i, (v,e_i')p_i \colon i \in [m]\}, \1 \\
E_4 & = & \{(v,s_i)s_i \colon i \in [n] \}, \1 \\
\end{array}
\]
we have $E(G') = E_1 \cup E_2 \cup E_3 \cup E_4$. An illustration of the construction of the graph $G'$ from the graph $G$ is given in Fig. \ref{MiF}.

\vspace{-6mm}
\begin{figure}[h!]
\begin{center} \bigskip
{\epsfxsize=6in \epsffile{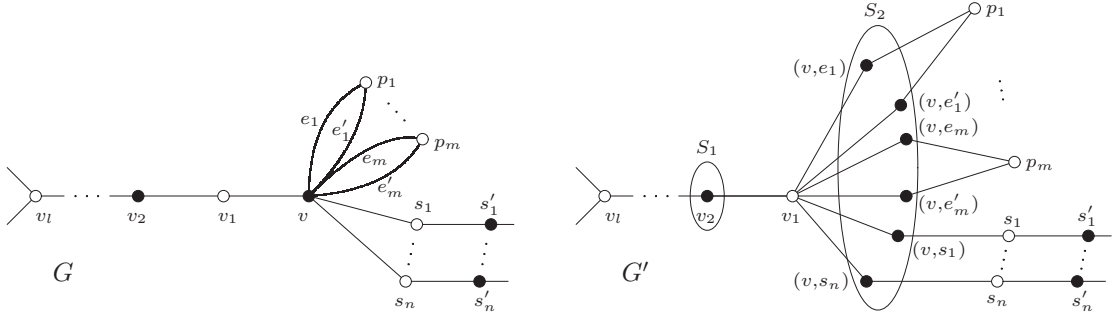}}

\vspace{-2mm}
 \caption{Graphs $G$ and $G'$ such that $G'({\cal P}(v_1))$ is isomorphic to $G$ } \label{MiF}
\end{center}
\end{figure}

We note that every new vertex added to $G$ when constructing $G'$ has degree~$2$ in $G'$. Further, the degrees of all vertices in $G$ different from $v$ remain unchanged in $G'$, except for the vertex $v_1$ whose degree changes from~$2$ to~$1 + 2m + n$. It follows from the fact that $G$ is a minimal $(2,2)$-dominated graph and from Theorem \ref{thm-property-22} that $G'$ is a minimal $(2,2)$-dominated graph, that is, from the fact that $G \in {\cal M}$ it follows that $G' \in {\cal M}$. Recall that by assumption, we have $|A_G^3| \ge |B_G^3|=k$. Since $B_{G'}^3 = B_G^3\setminus \{v\}$ and $A_{G'}^3= A_G^3\cup \{v_1\}$, we therefore have
\[
\min\{|A_{G'}^3|, |B_{G'}^3|\} = \min\{|A_{G}^3|+1, |B_{G}^3|-1\} = |B_{G}^3|-1 = k-1<k.
\]

Applying the induction hypothesis to the graph $G'$, we have that $G' \in {\cal F}$. Finally, if ${\cal P}(v_1)=\{S_1,S_2\}$ is a partition of $N_{G'}(v_1)$, where $S_1=\{v_2\}$ and $S_2=  N_{G'}(v_1)\setminus \{v_2\}= V^*$, then the $\cal P$-contraction $G'({\cal P}(v_1))$ belongs to the family $\cal F$. Consequently, the graph $G$ belongs to $\cal F$ as $G$ is isomorphic to $G'({\cal P}(v_1))$. This completes the proof of Theorem~\ref{char:M}.
\end{proof}

\section{Acknowledgements}

The authors wish to thank Dieter Rautenbach whose valuable comments improved
the clarity of the paper.

\end{document}